\documentclass[12pt,oneside,reqno]{amsart}
\usepackage{graphicx}
\usepackage{mathrsfs}
\usepackage{stmaryrd}
\usepackage{amsfonts}
\usepackage{cite}
\usepackage{enumerate,amsmath,amssymb,amsthm}
\usepackage{booktabs} 
\usepackage{diagbox} 

\newcommand{\dif}{\mathrm{d}}

\newcommand{\be}{\begin{eqnarray}}
	\newcommand{\ee}{\end{eqnarray}}
\newcommand{\ce}{\begin{eqnarray*}}
	\newcommand{\de}{\end{eqnarray*}}
\newtheorem{theorem}{Theorem}[section]
\newtheorem{lemma}[theorem]{Lemma}
\newtheorem{remark}[theorem]{Remark}
\newtheorem{definition}[theorem]{Definition}
\newtheorem{proposition}[theorem]{Proposition}
\newtheorem{Examples}[theorem]{Examples}
\newtheorem{corollary}[theorem]{Corollary}
\newtheorem{condition}[theorem]{Condition}
\def\e{\varepsilon}
\def\t{\theta}
\def\a{\alpha}

\def\b{\beta}

\def\p{\partial}

\def\s{\sigma}

\def\l{\lambda}

\def\[{{\Big[}}
\def\]{{\Big]}}
\def\<{{\langle}}
\def\>{{\rangle}}
\def\({{\Big(}}
\def\){{\Big)}}

\def\tr{{\rm tr}}

\def\no{\nonumber}
\def\bt{\begin{theorem}}
	\def\et{\end{theorem}}
\def\bl{\begin{lemma}}
	\def\el{\end{lemma}}
\def\br{\begin{remark}}
	\def\er{\end{remark}}
\def\bx{\begin{Examples}}
	\def\ex{\end{Examples}}
\def\bd{\begin{definition}}
	\def\ed{\end{definition}}
\def\bp{\begin{proposition}}
	\def\ep{\end{proposition}}
\def\bc{\begin{corollary}}
	\def\ec{\end{corollary}}
\def\bco{\begin{condition}}
	\def\eco{\end{condition}}

\def\mE{{\mathbb E}}

\def\mN{{\mathbb N}}

\def\mP{{\mathbb P}}

\def\mR{{\mathbb R}}

\def\mU{{\mathbb U}}
\def\mV{{\mathbb V}}

\def\sB{{\mathscr B}}

\def\sF{{\mathscr F}}

\def\sL{{\mathscr L}}

\def\geq{\geqslant}
\def\leq{\leqslant}

\begin{document}
	
\allowdisplaybreaks
\title{Large deviations for generalized backward stochastic differential equations}
	
\author{Yawen Liu and Huijie Qiao$^*$}

	
	
\thanks{This work was supported by NSF of China (No.12071071) and the Jiangsu Provincial Scientific Research Center of Applied Mathematics (No. BK20233002).}

\thanks{$*$ Corresponding author: Huijie Qiao, hjqiaogean@seu.edu.cn}
	
	
\date{}
	
\dedicatory{
Department of Mathematics, Southeast University,\\
Nanjing, Jiangsu 211189, P.R.China}
	
\begin{abstract}
This work concerns generalized backward stochastic differential equations, which are coupled with a family of reflecting diffusion processes. First of all, we establish the large deviation principle for forward stochastic differential equations with reflecting boundaries under weak monotonicity conditions. Then based on the obtained result and the contraction principle, the large deviation principle for the generalized backward stochastic differential equations is proved. As a by-product, we obtain a limit result about parabolic partial differential equations with the nonlinear Neumann boundary conditions. 	
\end{abstract}
	
\maketitle \rm
	
\section{Introduction}

Backward stochastic differential equations (BSDEs for short), initially introduced by Pardoux and Peng \cite{pp}, provide probabilistic formulas for the solutions of semilinear partial differential equations of both parabolic and elliptic types. Subsequent to their introduction, BSDEs have garnered considerable attention within the fields of mathematical finance, stochastic optimal control, partial differential equations and so on (cf. \cite{de, kp, lt, mz, pr1, pr2, pz, sp, zj}, and the references cited therein). 

The asymptotic theory of large deviation principles (LDPs for short) aims to make exponential asymptotic estimates for the probability of rare events, and LDPs have been a fundamental aspect of applied probability research. LDPs for BSDEs have been established for some cases (cf. \cite{ee,km,mz,nd,rs,sw}). Let us mention some works related with our results. In \cite{ee} Essaky studied the LDP for the following BSDE on $\mR^k$: for every $0\leq s\leq t\leq T$,
\be
\left\{\begin{array}{ll}
\bar Y_t^{s,x,\varepsilon}=h(\bar X_T^{s,x,\varepsilon})+\int_t^T f(r,\bar X_r^{s,x,\varepsilon},\bar Y_r^{s,x,\varepsilon}, \bar Z_r^{s,x,\varepsilon})\dif r-\int_t^T\bar Z_r^{s,x,\varepsilon}\dif W_r-\int_t^T\bar U_r^{s,x,\varepsilon}\dif r,\\
(\bar Y_t^{s,x,\varepsilon},\bar U_t^{s,x,\varepsilon})\in\p\Sigma, \quad \mE\int_0^T\Sigma(\bar Y_t^{s,x,\varepsilon})\dif t<\infty,
\end{array}
\right.
\label{barymbsde}
\ee
where $W=\{W_t,t\geq 0\}$ is a $m$-dimensional Brownian motion defined on a complete probability space $(\Omega, \sF, \mP)$, $\p\Sigma$ denotes the subdifferential operator of the proper lower semicontinuous convex function $\Sigma$, $h: \bar{\varTheta} \mapsto \mR^k$, $f: [0,T]\times\bar{\varTheta} \times\mR^k\times\mR^{k\times m}\mapsto \mR^k$ are Borel measurable, $\bar X^{s,x,\varepsilon}$ is the solution of the following stochastic differential equation (SDE for short) with the reflecting boundary on $\mR^d$:
\be
\left\{\begin{array}{ll}
\bar X_t^{s,x,\varepsilon}=x+\int_s^t \bar b(\bar X_r^{s,x,\varepsilon})\dif r+\sqrt{\varepsilon}\int_s^t\bar \sigma(\bar X_r^{s,x,\varepsilon})\dif W_r+\int_s^t\nabla\phi(\bar X_r^{s,x,\varepsilon})\dif \bar K_r^{s,x,\varepsilon},\\
\bar K_t^{s,x,\varepsilon}=\int_s^t\mathbb{I}_{\{\bar X_r^{s,x,\varepsilon}\in\partial\varTheta\}}\dif \bar K_r^{s,x,\varepsilon},\quad \bar K_{.}^{s,x,\varepsilon}\ is\ increasing,
\end{array}
\right.
\label{reflsde}
\ee
$x\in \bar{\varTheta}$, $\varepsilon>0$ is a small parameter, $\bar b: \bar{\varTheta} \mapsto \mR^d$, $\bar\sigma: \bar{\varTheta} \mapsto \mR^{d\times m}$ are Borel measurable and $\varTheta$ is an open connected bounded convex subset of $\mathbb{R}^d$ such that for a function $\phi\in C_b^2(\mR^d)$ (See Subsection \ref{nota}), $\varTheta=\{\phi>0\}$, $\partial\varTheta=\{\phi=0\}$, and $|\nabla\phi(x)|=1$, $x\in\partial\varTheta$. Note that at any boundary point $x\in\partial\varTheta$, $\nabla\phi(x)$ is a unit normal vector to the boundary, pointing towards the interior of $\varTheta$. Later, N'Zi and Dakaou \cite{nd} deleted the reflecting boundary in Eq.(\ref{reflsde}), replaced $\p\Sigma$ by any multivalued maximal monotone operator $A$ and also obtained the LDP result.

In this paper, we consider the following generalized backward coupled stochastic differential equation on $\mR^k$: for every $0\leq s\leq t\leq T$,
\be
Y_t^{s,x,\varepsilon}&=&h(X_T^{s,x,\varepsilon})+\int_t^Tf(r,X_r^{s,x,\varepsilon},Y_r^{s,x,\varepsilon},Z_r^{s,x,\varepsilon})\dif r+\int_t^Tg(r,X_r^{s,x,\varepsilon},Y_r^{s,x,\varepsilon})\dif K_r^{s,x,\varepsilon}\no\\
&&-\int_t^TZ_r^{s,x,\varepsilon}\dif W_r,
\label{ygbsde}
\ee
where $h, f$ are defined as above, $g: [0,T]\times\bar{\varTheta} \times\mR^k \mapsto \mR^k$ is Borel measurable, $X^{s,x,\varepsilon}$ is the solution of the following SDE with the reflecting boundary on $\mR^d$:
\be
\left\{\begin{array}{ll}
X_t^{s,x,\varepsilon}=x+\int_s^t b(r,X_r^{s,x,\varepsilon})\dif r+\sqrt{\varepsilon}\int_s^t\sigma(r,X_r^{s,x,\varepsilon})\dif W_r+\int_s^t\nabla\phi(X_r^{s,x,\varepsilon})\dif K_r^{s,x,\varepsilon},\\
K_t^{s,x,\varepsilon}=\int_s^t\mathbb{I}_{\{X_r^{s,x,\varepsilon}\in\partial\varTheta\}}\dif K_r^{s,x,\varepsilon},\quad K_{.}^{s,x,\varepsilon}\ is\ increasing,
\end{array}
\right.
\label{reflsdet}
\ee
and $b: [0,T]\times\bar{\varTheta} \mapsto \mR^d$, $\sigma: [0,T]\times\bar{\varTheta} \mapsto \mR^{d\times m}$ are Borel measurable. Under the suitable assumption, Eq.(\ref{ygbsde}) has a unique solution $(Y^{s,x,\varepsilon},Z^{s,x,\varepsilon})$. Moreover, $u^\varepsilon(s,x):=Y_s^{s,x,\varepsilon}$ is a unique viscosity solution of the following partial differential equation (PDE for short):
\be\left\{\begin{array}{l}
\frac{\p u_i^\e(s,x)}{\p s}+\sL_s^\e u_i^\e(s,x)+f_i(s, x, u^\e(s,x), (\triangledown u_i^\e \s)(s,x) )=0,\\
 ~(s,x)\in[0,T]\times\bar{\varTheta}, \quad i=1, \cdots, k\\ 
\frac{\p u_i^\e(s,x)}{\p n}+g_i(s,x,u^\e(s,x))=0, \\
~(s,x)\in[0,T]\times\p \bar{\varTheta}, \quad i=1, \cdots, k\\
u^\e(T,x)=h(x), \quad x\in\bar{\varTheta},
\end{array}
\right.
\label{pde1}
\ee
where 
\ce
&&(\sL_s^\e u_i^\e)(s,x)=b_j(s,x)\frac{\p u_i^\e(s,x)}{\p x_j}+\frac{\e}{2}(\s\s^*)_{jl}(s,x) \frac{\p^2 u_i^\e(s,x)}{\p x_j\p x_l}, \\
&&\frac{\p u_i^\e(s,x)}{\p n}=-\<\triangledown \phi(x), \triangledown u_i^\e(s,x)\>,
\de
and we use the convention that repeated indices imply summation. Note that Eq.(\ref{pde1}) is a parabolic PDE with the nonlinear Neumann boundary condition. And as usual it is not easy to solve this type of parabolic PDEs and let alone study the properties of their solutions. Thus it is an optional approach to investigate properties for the solutions of Eq.(\ref{ygbsde}) so as to obtain the properties of solutions to the system (\ref{pde1}). In \cite{lt, pr1, pz}, authors have established the well-posedness for the solutions of Eq.(\ref{ygbsde}). Therefore, in this paper, our purpose is to study the asymptotic behavior of the family $(Y^{s,x,\varepsilon})_{\varepsilon>0}$ as $\varepsilon$ goes to 0. We firstly establish the LDP for forward SDEs with reflecting boundaries under weak monotonicity conditions. And our result can cover \cite[Theorem 2.1]{ee}. Then based on the obtained result and the contraction principle, the large deviation principle for Eq.(\ref{ygbsde}) is proved. Moreover, we obtain a limit result about parabolic PDEs with the nonlinear Neumann boundary conditions (See Remark \ref{limi}). 	

Comparing Eq.(\ref{reflsdet}) and Eq.(\ref{ygbsde}) with Eq.(\ref{reflsde}) and Eq.(\ref{barymbsde}), we find that Eq.(\ref{reflsdet}) is more general than Eq.(\ref{reflsde}), and there is the term $\int_t^Tg(r,X_r^{s,x,\varepsilon},Y_r^{s,x,\varepsilon})\dif K_r^{s,x,\varepsilon}$ in Eq.(\ref{ygbsde}). Since $\dif K_r^{s,x,\varepsilon}$ is not absolutely continuous with respect to the Lebesgue measure, it is difficult to estimate this term straightly. Here in order to obtain the LDP result we use some new techniques to overcome this difficulty.

The rest of this paper is organized as follows. In Section \ref{pre} we introduce some notations, concepts and some related results about LDPs. In Section \ref{sdeLDP} the LDP for forward SDEs with reflecting boundaries is presented. Section \ref{bsdeLDP} focuses on the LDP for generalized BSDEs.
	
The following convention will be used throughout the paper: $C$ with or without indices will denote different positive constants whose values may change from one place to another.
	
\section{Preliminaries}\label{pre}

In the section, we introduce notations and concepts.

\subsection{Notations}\label{nota}

In this subsection, we present notes and notations used in the sequel.
	
Let $|\cdot|$ and $\|\cdot\|$ be norms of vectors and matrices, respectively. Let $\langle\cdot,\cdot\rangle$ be the inner product in $\mathbb{R}^d$. 

$C^2(\mR^d)$ stands for the space of continuous functions on $\mR^d$ which have continuous partial derivatives of order up to $2$, and $C_b^2(\mR^d)$ stands for the subspace of $C^2(\mR^d)$, consisting of functions whose derivatives up to order 2 are bounded.

$C([0,T], \mathbb{R}^d)$ is the set of continuous functions on $[0,T]$ with values in $\mathbb{R}^d$. We equip $C([0,T], \mathbb{R}^d)$ with the uniform convergence topology. For any $\varrho\in C([0,T], \mathbb{R}^d)$, $\|\varrho\|_\infty:= \sup \limits_{t\in [s,T]}|\varrho_t|$. $V_0([0,T], \mR^d)$ denotes the space of functions $\rho: [0,T]\to \mR^d$ with bounded variation and $\rho_0=0$. For $\rho\in V_0([0,T], \mR^d)$ and $s\in [0,T]$, we shall use $|\rho|_{0}^{s}$ to denote the variation of $\rho$ on $[0,s]$ and write $|\rho|_{TV}:=|\rho|_{0}^{T}$.

\subsection{Large deviation principles}
	
In this subsection, we define large deviation principles and give some related results.

Let $\mU$ be a Polish space. Let $\{X^{\e}, \e>0\}$ be a family of $\mU$-valued random variables defined on $(\Omega, \mathscr{F}, \{\mathscr{F}_t\}_{t\in[0,T]}, \mP)$.

\bd
$(i)$ A function $\Lambda: \mU\rightarrow[0,+\infty]$ is called a rate function on $\mU$, if for all $M \geq 0$, $\{x\in\mU: \Lambda(x)\leq M \}$ is a closed subset of $\mU$.
	
$(ii)$ A rate function $\Lambda$ is called a good rate function, if for all $M\geq 0$, $\{x\in\mU: \Lambda(x)\leq M \}$ is a compact subset of $\mU$.
\ed
	
\bd
$\{X^{\e}, \e>0\}$ is said to satisfy the large deviation principle with the good rate function $\Lambda$ if for any subset $B\in\sB(\mU)$
$$
-\inf\limits_{x \in B^\circ} \Lambda(x) \leq \liminf\limits_{\varepsilon \to 0} \varepsilon \ln\(\mP(X^\varepsilon \in B^\circ)\) \leq \limsup\limits_{\varepsilon \to 0} \varepsilon \ln\(\mP(X^\varepsilon \in \bar B)\) \leq -\inf\limits_{x\in \bar{B}} \Lambda(x),
$$
where $B^\circ, \bar{B}$ are the interior and closure of $B$ in $\mU$, respectively.
\ed

The following contraction principle is from \cite[Theorem 2.4, P.5]{va} (See \cite[Theorem 4.2.1, P.126]{dz} for the case $F^\e=F$).

\bt\label{yasuo}
Let $\mU$ and $\mV$ be two Polish spaces and $F^\e: \mU\rightarrow\mV$ be a continuous function. Assume that 
$\lim\limits_{\e\rightarrow0}F^\e=F$ exists uniformly over compact subsets of $\mU$. Consider a good rate function $\Lambda: \mU \rightarrow [0,\infty]$.

$(i)$ For each $y\in\mV$, define
$$
S(y) \triangleq \inf\{\Lambda(x): x \in\mU,\quad y=F(x) \}.
$$
Then $S$ is a good rate function on $\mV$, where as usual the infimum over the empty set is taken as $\infty$.

$(ii)$ If $\Lambda$ controls the large deviation principle associated with a family of probability measures $\{\mu^\varepsilon\}$ on $\mU$, then $S$ controls the large deviation principle associated with the family of probability measures $\{\mu^\varepsilon\circ (F^\e)^{-1}\}$ on $\mV$.
\et

\section{The LDP and convergence for forward SDEs with reflecting boundaries}\label{sdeLDP}

In this section, we prove the LDP and convergence for forward SDEs with reflecting boundaries.

First of all, we mention that the assumption on $\varTheta$ implies 
\be
2\<x'-x,\nabla\phi(x)\>+\a|x-x'|^2\geq 0, \quad x\in\p\varTheta, x'\in\bar \varTheta
\label{fielcond}
\ee
where $\a>0$ is a constant. Besides, we define the following function associated with $\bar\varTheta$
$$
I_{\bar\varTheta}(x):= 
\begin{cases}0, & \text { if } x \in \bar\varTheta, \\ 
+\infty, & \text { if } x \notin \bar\varTheta.
\end{cases}
$$
The subdifferential operator of $I_{\bar\varTheta}$ is given by
$$
\begin{aligned}
\partial I_{\bar\varTheta}(x) & :=\left\{y \in \mathbb{R}^d:\langle y, z-x\rangle \leq 0, \forall z \in \bar\varTheta\right\} \\
& = \begin{cases}\emptyset, & \text { if } x \notin \bar\varTheta, \\
\{0\}, & \text { if } x \in \operatorname{Int}(\bar\varTheta), \\
\Lambda_x, & \text { if } x \in \partial \bar\varTheta,\end{cases}
\end{aligned}
$$
where $\Lambda_x$ is the exterior normal cone at $x$. By simple deduction, we know that $\partial I_{\bar\varTheta}$ is a maximal monotone operator.
	
Let $\{\sF_{t}\}_{t\geq0}$ be the augmented natural filtration of $W$. On the complete filtered probability space $(\Omega, \sF, \mP, \{\sF_{t}\}_{t\geq0})$ we consider Eq.(\ref{reflsdet}), i.e.
\be
\left\{\begin{array}{ll}
\dif X_t^{s,x,\varepsilon}\in \p I_{\bar\varTheta}(X_t^{s,x,\varepsilon})\dif t+b(t,X_t^{s,x,\varepsilon})\dif t+\sqrt{\varepsilon}\sigma(t,X_t^{s,x,\varepsilon})\dif W_t,\\
X_s^{s,x,\varepsilon}=x,
\end{array}
\right.
\label{sdef1}
\ee
where $x\in \bar{\varTheta}$ and $b: \bar{\varTheta} \mapsto \mR^d$, $\sigma: \bar{\varTheta} \mapsto \mR^{d\times m}$ are Borel measurable.

\subsection{The LDP for forward SDEs with reflecting boundaries}

In this subsection, we prove the LDP for Eq.(\ref{sdef1}) under weak monotonicity conditions.

First of all, for Eq.(\ref{sdef1}), we assume:
\begin{enumerate}[(${\bf H}^1_{b,\s}$)]
\item $(i)$ $b(t,x), \s(t,x)$ are continuous in $x$, and there exists $\t\in(0,1)$ such that for any $t\in[0,T], x, x'\in\mR^d, |x-x'|\leq \t$
\ce
2\<x-x', b(t,x)-b(t,x')\>+\|\s(t,x)-\s(t,x')\|^2\leq \kappa_1(t)\gamma_1(|x-x'|^2),
\de
where $\kappa_1$ is a positive integrable function on $[0,T]$ and $\gamma_1: [0,1)\rightarrow \mR_+$ is an increasing continuous function satisfying 
\ce
\gamma_1(0)=0, \quad \int_{0^+}\frac{\dif x}{\gamma_1(x)}=\infty, 
\de
and for any constant $0<\zeta_1\leq 1$
$$
\sup\limits_{x\in [0,\t]}\frac{\zeta_1 \gamma_1(x)}{\gamma_1(\zeta_1 x)}<\infty.
$$
$(ii)$ There exist two constants $c_1, c_2>0$ such that for any $t\in[0,T], x\in\mR^d$
\ce
2\<b(t,x), x\>+c_1\|\s(t,x)\|^2+c_2\frac{|\<\s(t,x), x\>|^2}{|x|^2}\leq \kappa_2(t)(1+\gamma_2(|x|^2)),
\de
where $\kappa_2$ is a positive integrable function on $[0,T]$ and $\gamma_2: [0, \infty)\rightarrow \mR_+$ is an increasing continuous function satisfying 
$$
\int_0^\infty\frac{\dif x}{\gamma_2(x)+1}=\infty,
$$
and for any constant $0<\zeta_2\leq 1$
$$
\sup\limits_{x\in [0,\infty)}\frac{\zeta_2 \gamma_2(x)}{\gamma_2(\zeta_2 x)}<\infty.
$$
$(iii)$
\ce
\int_0^T(|b(t,x)|+\|\s(t,x)\|^2)\dif t<\infty.
\de
\end{enumerate}
\begin{enumerate}[(${\bf H}^2_{\s}$)]
\item There exists a constant $\iota>0$ such that for any $t\in[0,T], x\in\mR^d$ and $h \in \mathbb{R}^d$
$$
\langle \s\s^*(t,x)h,h \rangle \geq \iota|h|^2.
$$
\end{enumerate}

Under $({\bf H}^1_{b,\s})$, it follows from the similar deduction to that in \cite[Theorem 2.8]{rwz} that there exists a unique pair of progressively measurable continuous processes $\{(X_t^{s,x,\varepsilon},K_t^{s,x,\varepsilon})\}_{0\leq s\leq t\leq T}$, with values in $\bar{\varTheta}\times\mR_+$ satisfying Eq.(\ref{sdef1}). That is,
\ce
\left\{\begin{array}{ll}
X_t^{s,x,\varepsilon}=x+\int_s^t b(r, X_r^{s,x,\varepsilon})\dif r+\sqrt{\varepsilon}\int_s^t\sigma(r, X_r^{s,x,\varepsilon})\dif W_r+\int_s^t\nabla\phi(X_r^{s,x,\varepsilon})\dif K_r^{s,x,\varepsilon},\\
K_t^{s,x,\varepsilon}=\int_s^t\mathbb{I}_{\{X_r^{s,x,\varepsilon}\in\partial\varTheta\}}\dif K_r^{s,x,\varepsilon},\quad K_{.}^{s,x,\varepsilon}\ is\ increasing.
\end{array}
\right.
\de

Next, we study the LDP for $X^{s,x,\varepsilon}$. In order to do this, we make preparations. Let $\Phi \in C([s,T], \mR^d)$, $\Psi \in C([s,T], \bar{\varTheta})$, $\rho \in V_0([s,T], \mR^d)$ such that
$$
\Psi_t=\Phi_t+\rho_t,\quad \rho_t=\int_s^t\nabla\phi(\Psi_r)\dif |\rho|^r_s,\quad |\rho|^t_s=\int_s^t\mathbb{I}_{\{\Psi_r\in\partial\varTheta\}}\dif |\rho|^r_s.
$$
Define the mapping $F: C([s,T], \mR^d)\rightarrow C([s,T], \bar{\varTheta})$ by $\Psi=F(\Phi)$. Then $F$ is well defined and continuous (\cite[Lemma 2.3]{sa}). 

Now, we state and prove the main result of this section.
		
\bt\label{fldp} 
Under the assumptions $({\bf H}^1_{b,\s})$ and $({\bf H}^2_{\s})$, $\{X^{s,x,\varepsilon}, \e>0\}$ satisfies a large deviation principle with the good rate function 
\ce
S(\Psi)=
\left\{\begin{array}{ll}
\frac{1}{2} \inf\limits_{\Phi \in F^{-1}(\Psi)} \int_s^T (\dot{\Phi}_t-b(t,\Psi_t))^* (\s\s^*)^{-1}(t,\Psi_t) (\dot{\Phi}_t-b(t,\Psi_t)) \dif t, \\
\qquad\quad\text{if $\Phi_t$ is absolutely continuous, $\Phi_s=x$ and $F^{-1}(\Psi) \neq \emptyset$};\\
+\infty, \quad \text{in the opposite case}. 
\end{array}
\right. 
\de
\et
\begin{proof}
First of all, we consider the following SDE on $\mR^d$:
\be
\xi^{s,x,\varepsilon}_t=x+\int_s^t b(r,\xi_r^{s,x,\varepsilon})\dif r+\sqrt{\varepsilon}\int_s^t\sigma(r,\xi_r^{s,x,\varepsilon})\dif W_r.
\label{sde}
\ee
Under $({\bf H}^1_{b,\s})$, Eq.(\ref{sde}) has a unique strong solution $\xi^{s,x,\varepsilon}$ with $\xi_s^{s,x,\varepsilon}=x$ (cf. \cite[Proposition 2.1]{wyzz}). Moreover, by \cite[Theorem 2.1]{wyzz}, we know that under $({\bf H}^2_{\s})$, $\{\xi^{s,x,\varepsilon}, \e>0\}$ satisfies the LDP with the good rate function
\ce
\Lambda(\Phi)=
\left\{\begin{array}{ll}
\frac{1}{2} \int_s^T (\dot{\Phi}_t-b(t,\Phi_t))^*(\s\s^*)^{-1}(t,\Phi_t)(\dot{\Phi}_t-b(t,\Phi_t))\dif t, \\
\qquad\quad\text{if $\Phi_t$ is absolutely continuous and $\Phi_s=x$};\\
+\infty, \quad \text{in the opposite case}.   
\end{array}
\right.
\de

Next, we notice that $X^{s,x,\varepsilon}=F(\xi^{s,x,\varepsilon})$ and $F$ is continuous. Therefore, from Theorem \ref{yasuo}, it follows that $\{X^{s,x,\varepsilon}, \e>0\}$ satisfies the LDP with the good rate function 
\ce
S(\Psi)
&=& \inf \{ \Lambda(\Phi):\Phi \in C([s,T], \mR^d),\quad \Psi=F(\Phi) \}\no\\
&=&\left\{\begin{array}{ll}
\frac{1}{2} \inf\limits_{\Phi \in F^{-1}(\Psi)} \int_s^T (\dot{\Phi}_t-b(t,\Psi_t))^* (\s\s^*)^{-1}(t,\Psi_t) (\dot{\Phi}_t-b(t,\Psi_t)) \dif t, \\
\qquad\quad\text{if $\Phi_t$ is absolutely continuous, $\Phi_s=x$ and $F^{-1}(\Psi) \neq \emptyset$};\\
+\infty, \quad \text{in the opposite case}. 
\end{array}
\right. 
\de
The proof is complete.
\end{proof}

\subsection{The convergence for forward SDEs with reflecting boundaries}

In this subsection, we give some convergence results for Eq.(\ref{sdef1}) under Lipschitz continuous conditions.

We assume:	
\begin{enumerate}[(${\bf H}^{1'}_{b,\s}$)]
\item There exists a constant $L_1>0$ such that for any $t\in[0,T]$ and $x,x'\in \bar{\varTheta}$
$$
|b(t, x)-b(t, x')|+\|\sigma(t,x)-\sigma(t,x')\|\leq L_1|x-x'|,
$$
and
$$
|b(t,x)|+\|\sigma(t,x)\|\leq L_1(1+|x|). 
$$
\end{enumerate}

\br
$({\bf H}^{1'}_{b,\s})$ implies $({\bf H}^{1}_{b,\s})$. Indeed, by $({\bf H}^{1'}_{b,\s})$, we take $\kappa_1(t)=2L_1+L^2_1, \gamma_1(x)=x, c_1=c_2=1, \kappa_2(t)=4L_1+4L_1^2, \gamma_2(x)=x$ , and obtain $({\bf H}^{1}_{b,\s})$.
\er

Based on the above remark and Theorem \ref{fldp}, we know that under $({\bf H}^{1'}_{b,\s})$ and $({\bf H}^2_{\s})$, the solution $X^{s,x,\varepsilon}$ of Eq.(\ref{sdef1}) satisfies the LDP.

In the following, in order to study the LDP about generalized BSDEs, we consider the ordinary differential equation with the reflecting boundary on $\mR^d$:
\be\left\{\begin{array}{l}
\dif \varphi_t^{s,x}\in \p I_{\bar\varTheta}\dif t+b(t,\varphi_t^{s,x})\dif t,\\
\varphi_s^{s,x}=x.
\end{array}
\right. 
\label{def1}
\ee
Under $({\bf H}^{1'}_{b,\s})$, the above equation (\ref{def1}) has a unique solution $(\varphi^{s,x}, K^{s,x})$ with values in $\bar{\varTheta}\times\mR_+$ satisfying
\ce
\left\{\begin{array}{ll}
\varphi_t^{s,x}=x+\int_s^t b(r, \varphi_r^{s,x})\dif r+\int_s^t\nabla\phi(\varphi_r^{s,x})\dif K_r^{s,x},\\
K_t^{s,x}=\int_s^t\mathbb{I}_{\{\varphi_r^{s,x}\in\partial\varTheta\}}\dif K_r^{s,x},\quad K_{.}^{s,x}\ is\ increasing.
\end{array}
\right. 
\de
We study the relationship between Eq.(\ref{sdef1}) and Eq.(\ref{def1}).
		
\bl\label{xshou}
Under $({\bf H}^{1'}_{b,\s})$, it holds that for $0<\e<1$
\ce
\mE\left[\sup\limits_{s \leq t \leq T}|X_t^{s,x,\varepsilon}-\varphi_t^{s,x}|^4|\sF_s\right]\leq C\varepsilon,
\de
where the constant $C>0$ is independent of $\e$.
\el
\begin{proof}		
First of all, for Eq.(\ref{sdef1}) and Eq.(\ref{def1}), by the It\^o formula, it holds that
\ce
\dif |X_t^{s,x,\varepsilon}-\varphi_t^{s,x}|^2&=& 2\big\langle X_t^{s,x,\varepsilon}-\varphi_t^{s,x} , b(t,X_t^{s,x,\varepsilon})-b(t,\varphi_t^{s,x}) \big\rangle \dif t\no\\
&&+2\big\langle X_t^{s,x,\varepsilon}-\varphi_t^{s,x}, \sqrt{\varepsilon}\sigma(t,X_t^{s,x,\varepsilon})\dif W_t \big\rangle \no\\
&&+2\big \langle X_t^{s,x,\varepsilon}-\varphi_t^{s,x} , \nabla\phi(X_t^{s,x,\varepsilon})\dif K_t^{s,x,\varepsilon}-\nabla\phi(\varphi_t^{s,x})\dif K_t^{s,x} \big \rangle \no\\
&&+\|\sqrt{\varepsilon}\s(t,X_t^{s,x,\varepsilon})\|^2\dif t.
\de
Besides, applying the It\^o formula to $e^{-\alpha\phi(X_t^{s,x,\varepsilon})}$ and $e^{-\alpha\phi(\varphi_t^{s,x})}$, we have 
\ce
\dif e^{-\alpha\phi(X_t^{s,x,\varepsilon})}&=& -\alpha e^{-\alpha\phi(X_t^{s,x,\varepsilon})} \Big[\<\nabla\phi(X_t^{s,x,\varepsilon}),\dif X_t^{s,x,\varepsilon}\>-\frac{\alpha}{2} |\sqrt{\varepsilon}\s^*(t,X_t^{s,x,\varepsilon})\nabla\phi(X_t^{s,x,\varepsilon})|^2\dif t\no\\
&&+\frac{1}{2}\tr\left(\sqrt{\varepsilon}\s^*(t,X_t^{s,x,\varepsilon})D^2\phi(X_t^{s,x,\varepsilon})\sqrt{\varepsilon}\s(t,X_t^{s,x,\varepsilon})\right)\dif t\Big],\\
\dif e^{-\alpha\phi(\varphi_t^{s,x})}&=& -\alpha e^{-\alpha\phi(\varphi_t^{s,x})}\<\nabla\phi(\varphi_t^{s,x}),\dif \varphi_t^{s,x}\>,
\de
which yields that
\ce
\dif e^{-\alpha\left(\phi(X_t^{s,x,\varepsilon})+\phi(\varphi_t^{s,x})\right)}
&=&e^{-\alpha\phi(\varphi_t^{s,x})} \dif e^{-\alpha\phi(X_t^{s,x,\varepsilon})}+e^{-\alpha\phi(X_t^{s,x,\varepsilon})} \dif e^{-\alpha\phi(\varphi_t^{s,x})}\\
&=& -\alpha e^{-\alpha\left(\phi(X_t^{s,x,\varepsilon})+\phi(\varphi_t^{s,x})\right)}\Big[\<\nabla\phi(X_t^{s,x,\varepsilon}),\dif X_t^{s,x,\varepsilon}\>+\<\nabla\phi(\varphi_t^{s,x}),\dif \varphi_t^{s,x}\>\no\\
&&-\frac{\alpha}{2} |\sqrt{\varepsilon}\s^*(t,X_t^{s,x,\varepsilon})\nabla\phi(X_t^{s,x,\varepsilon})|^2\dif t\no\\
&&+\frac{1}{2} \tr\left(\sqrt{\varepsilon}\s^*(t,X_t^{s,x,\varepsilon})D^2\phi(X_t^{s,x,\varepsilon})\sqrt{\varepsilon}\s(t,X_t^{s,x,\varepsilon})\right)\dif t\Big].
\de

Collecting the above deduction, by the It\^o formula, $({\bf H}^{1'}_{b,\s})$ and (\ref{fielcond}), we obtain that
\ce
&&\dif e^{-\alpha\left(\phi(X_t^{s,x,\varepsilon})+\phi(\varphi_t^{s,x})\right)}|X_t^{s,x,\varepsilon}-\varphi_t^{s,x}|^2\no\\
&=& e^{-\alpha\left(\phi(X_t^{s,x,\varepsilon})+\phi(\varphi_t^{s,x})\right)} \dif |X_t^{s,x,\varepsilon}-\varphi_t^{s,x}|^2 + |X_t^{s,x,\varepsilon}-\varphi_t^{s,x}|^2 \dif e^{-\alpha\left(\phi
(X_t^{s,x,\varepsilon})+\phi(\varphi_t^{s,x})\right)}\no\\
&&+ \dif |X_t^{s,x,\varepsilon}-\varphi_t^{s,x}|^2 \dif e^{-\alpha\left(\phi(X_t^{s,x,\varepsilon})+\phi(\varphi_t^{s,x})\right)}\no\\
&=&e^{-\alpha\left(\phi(X_t^{s,x,\varepsilon})+\phi(\varphi_t^{s,x})\right)} 2\big\langle X_t^{s,x,\varepsilon}-\varphi_t^{s,x}, b(t,X_t^{s,x,\varepsilon})-b(t,\varphi_t^{s,x}) \big\rangle \dif t\no\\
&&+ 2e^{-\alpha\left(\phi(X_t^{s,x,\varepsilon})+\phi(\varphi_t^{s,x})\right)}\big\langle X_t^{s,x,\varepsilon}-\varphi_t^{s,x}, \sqrt{\varepsilon}\sigma(t,X_t^{s,x,\varepsilon})\dif W_t \big\rangle \no\\
&&+2e^{-\alpha\left(\phi(X_t^{s,x,\varepsilon})+\phi(\varphi_t^{s,x})\right)}\big \langle X_t^{s,x,\varepsilon}-\varphi_t^{s,x}, \nabla\phi(X_t^{s,x,\varepsilon})\dif K_t^{s,x,\varepsilon}-\nabla\phi(\varphi_t^{s,x})\dif K_t^{s,x} \big \rangle \no\\
&&+e^{-\alpha\left(\phi(X_t^{s,x,\varepsilon})+\phi(\varphi_t^{s,x})\right)}\|\sqrt{\varepsilon}\s(t,X_t^{s,x,\varepsilon})\|^2\dif t\no\\
&&-\alpha e^{-\alpha\left(\phi(X_t^{s,x,\varepsilon})+\phi(\varphi_t^{s,x})\right)}|X_t^{s,x,\varepsilon}-\varphi_t^{s,x}|^2\Big[\<\nabla\phi(X_t^{s,x,\varepsilon}),\dif X_t^{s,x,\varepsilon}\>+\<\nabla\phi(\varphi_t^{s,x}),\dif \varphi_t^{s,x}\>\no\\
&&-\frac{\alpha}{2}|\sqrt{\varepsilon}\s^*(t,X_t^{s,x,\varepsilon})\nabla\phi(X_t^{s,x,\varepsilon})|^2\dif t+\frac{1}{2} \tr\left(\sqrt{\varepsilon}\s^*(t,X_t^{s,x,\varepsilon})D^2\phi(X_t^{s,x,\varepsilon})\sqrt{\varepsilon}\s(t,X_t^{s,x,\varepsilon})\right)\dif t\Big]\no\\
&&-2\alpha\e e^{-\alpha\left(\phi(X_t^{s,x,\varepsilon})+\phi(\varphi_t^{s,x})\right)}\big\langle X_t^{s,x,\varepsilon}-\varphi_t^{s,x},(\s\s^*)(t,X_t^{s,x,\varepsilon})\nabla\phi(X_t^{s,x,\varepsilon})\big\rangle\dif t\\
&\leq&2L_1 e^{-\alpha\left(\phi(X_t^{s,x,\varepsilon})+\phi(\varphi_t^{s,x})\right)}|X_t^{s,x,\varepsilon}-\varphi_t^{s,x}|^2\dif t\\
&&+2e^{-\alpha\left(\phi(X_t^{s,x,\varepsilon})+\phi(\varphi_t^{s,x})\right)}\big\langle X_t^{s,x,\varepsilon}-\varphi_t^{s,x}, \sqrt{\varepsilon}\sigma(t,X_t^{s,x,\varepsilon})\dif W_t \big\rangle \no\\
&&+C\e e^{-\alpha\left(\phi(X_t^{s,x,\varepsilon})+\phi(\varphi_t^{s,x})\right)}\dif t+C\alpha e^{-\alpha\left(\phi(X_t^{s,x,\varepsilon})+\phi(\varphi_t^{s,x})\right)}|X_t^{s,x,\varepsilon}-\varphi_t^{s,x}|^2\dif t\\
&&+2\alpha C\varepsilon e^{-\alpha\left(\phi(X_t^{s,x,\varepsilon})+\phi(\varphi_t^{s,x})\right)}|X_t^{s,x,\varepsilon}-\varphi_t^{s,x}|\dif t\\
&\leq&2L_1 |X_t^{s,x,\varepsilon}-\varphi_t^{s,x}|^2\dif t+2\sqrt{\varepsilon}e^{-\alpha\left(\phi(X_t^{s,x,\varepsilon})+\phi(\varphi_t^{s,x})\right)}\big\langle X_t^{s,x,\varepsilon}-\varphi_t^{s,x},\sigma(t,X_t^{s,x,\varepsilon})\dif W_t \big\rangle \no\\
&&+C\e\dif t+C\alpha|X_t^{s,x,\varepsilon}-\varphi_t^{s,x}|^2\dif t+2\alpha C\varepsilon |X_t^{s,x,\varepsilon}-\varphi_t^{s,x}|\dif t.
\de
Then the Burkholder-Davis-Gundy inequality and the H\"older inequality imply that
\ce
&&\mE\left[\sup\limits_{s \leq t \leq T}|X_t^{s,x,\varepsilon}-\varphi_t^{s,x}|^4|\sF_s\right]\\
&\leq& 20L^2_1TC\mE\left[\int_s^T|X_r^{s,x,\varepsilon}-\varphi_r^{s,x}|^4\dif r|\sF_s\right]\\
&&+C\e\mE\left[\int_s^T|\sigma^*(r,X_t^{s,x,\varepsilon})(X_r^{s,x,\varepsilon}-\varphi_r^{s,x})|^2\dif r|\sF_s\right]\\
&&+5C\e^2T+5C\alpha^2T\mE\left[\int_s^T|X_r^{s,x,\varepsilon}-\varphi_r^{s,x}|^4\dif r|\sF_s\right]\\
&&+20\alpha^2 C\varepsilon^2T\mE\left[\int_s^T|X_r^{s,x,\varepsilon}-\varphi_r^{s,x}|^2\dif r|\sF_s\right]\\
&\leq& C\e+C\int_s^T\mE\left[\sup\limits_{s \leq u \leq r}|X_u^{s,x,\varepsilon}-\varphi_u^{s,x}|^4|\sF_s\right]\dif r,
\de
where we use the fact that $|X_r^{s,x,\varepsilon}|\leq C$ ($\varTheta$ is bounded and $X_r^{s,x,\varepsilon}\in \bar{\varTheta}$ for any $r\in[s,T]$) and $|X_r^{s,x,\varepsilon}-\varphi_r^{s,x}|^2\leq 1+|X_r^{s,x,\varepsilon}-\varphi_r^{s,x}|^4$. Finally, the Gronwall inequality yields the required estimate. The proof is complete.
\end{proof}

Based on the above estimate, we can obtain the following result which is used in the sequel.

\bl\label{kshou}
Under $({\bf H}^{1'}_{b,\s})$, it holds that for $0<\e<1$
\ce
\mE\left[\sup \limits_{s \leq t \leq T}|K_t^{s,x,\varepsilon}-K_t^{s,x}|^4|\sF_s\right]\leq C\varepsilon,
\de
where the constant $C>0$ is independent of $\e$.
\el
\begin{proof}
First of all, by the It\^o formula, it holds that for $s\leq t\leq T$
\ce
K_t^{s,x,\varepsilon}&=&\phi(X_t^{s,x,\varepsilon})-\phi(x)-\int_s^t\<\nabla\phi(X_r^{s,x,\varepsilon}),b(r,X_r^{s,x,\varepsilon})\>\dif r\\
&&-\frac{\e}{2}\int_{s}^{t} \tr\(\s^*(r,X_r^{s,x,\varepsilon})D^2\phi(X_r^{s,x,\varepsilon})\s(r,X_r^{s,x,\varepsilon})\)\dif r\\
&&-\sqrt{\varepsilon}\int_{s}^{t}\<\nabla\phi(X_r^{s,x,\varepsilon}),\sigma(r,X_r^{s,x,\varepsilon})\dif W_r\>,\\
K_t^{s,x}&=&\phi(\varphi_t^{s,x})-\phi(x)-\int_{s}^{t}\<\nabla\phi(\varphi_r^{s,x}),b(r,\varphi_r^{s,x})\>\dif r.
\de
Then $({\bf H}^{1'}_{b,\s})$ and the H\"older inequality imply that
\ce
|K_t^{s,x,\varepsilon}-K_t^{s,x}|^4&\leq& 4|\phi(X_t^{s,x,\varepsilon})-\phi(\varphi_t^{s,x})|^4\\
&&+4\left|\int_s^t\left(\<\nabla\phi(X_r^{s,x,\varepsilon}),b(r,X_r^{s,x,\varepsilon})\>-\<\nabla\phi(\varphi_r^{s,x}),b(r,\varphi_r^{s,x})\>\right)\dif r\right|^4\\
&&+4\left(\frac{\e}{2}\right)^4\left|\int_{s}^{t} \tr\(\s^*(r,X_r^{s,x,\varepsilon})D^2\phi(X_r^{s,x,\varepsilon})\s(r,X_r^{s,x,\varepsilon})\)\dif r\right|^4\\
&&+4\e^2\left|\int_{s}^{t}\<\nabla\phi(X_r^{s,x,\varepsilon}),\sigma(r,X_r^{s,x,\varepsilon})\dif W_r\>\right|^4\\
&\leq& 4C|X_t^{s,x,\varepsilon}-\varphi_t^{s,x}|^4+4CT^3\int_s^t|X_r^{s,x,\varepsilon}-\varphi_r^{s,x}|^4\dif r+C\e^4\\
&&+4\e^2\left|\int_{s}^{t}\<\nabla\phi(X_r^{s,x,\varepsilon}),\sigma(r,X_r^{s,x,\varepsilon})\dif W_r\>\right|^4,
\de
which together with the Burkholder-Davis-Gundy inequality and Lemma \ref{xshou} yields that
\ce
&&\mE\left[\sup \limits_{s \leq t \leq T}|K_t^{s,x,\varepsilon}-K_t^{s,x}|^4|\sF_s\right]\\
&\leq& 4C\mE\left[\sup \limits_{s \leq t \leq T}|X_t^{s,x,\varepsilon}-\varphi_t^{s,x}|^4|\sF_s\right]+4CT^3\int_s^T\mE[|X_r^{s,x,\varepsilon}-\varphi_r^{s,x}|^4|\sF_s]\dif r\\
&&+C\e^4+4\e^2 C\mE\left[\left(\int_{s}^{T}|\sigma^*(r,X_r^{s,x,\varepsilon})\nabla\phi(X_r^{s,x,\varepsilon})|^2\dif r\right)^2\bigg{|}\sF_s\right]\\
&\leq&C\e.
\de
The proof is complete.
\end{proof}

In the sequel, we also need the following estimate. 

\bl\label{kshou}
Assume that $({\bf H}^{1'}_{b,\s})$ holds. For all $\varepsilon \in (0,1)$ and $p\geq 1$, there exists a constant $C>0$ such that for all $x\in\bar{\varTheta}$,
\be
\mE\left[\sup\limits_{s \leq t \leq T}(K^{s,x,\varepsilon}_t)^p|\sF_s\right] \leq C, 
\label{k1}
\ee
and
\be
\sup\limits_{s \leq t \leq T}(K^{s,x}_t)^p\leq C.
\label{k2}
\ee
Moreover, for all $\b>0$ and $t\in[s,T]$, there exists a constant $C>0$ such that for all $x\in\bar{\varTheta}$,
\be
\mE [e^{\b K^{s,x,\varepsilon}_t}|\sF_s]\leq C.
\label{kexpo}
\ee
\el
\begin{proof}
The It\^o formula yields that
\ce
K_t^{s,x,\varepsilon}&=&\phi(X_t^{s,x,\varepsilon})-\phi(x)-\int_s^t\<\nabla\phi(X_r^{s,x,\varepsilon}),b(r,X_r^{s,x,\varepsilon})\>\dif r\\
&&-\frac{\e}{2}\int_{s}^{t} \tr\(\s^*(r,X_r^{s,x,\varepsilon})D^2\phi(X_r^{s,x,\varepsilon})\s(r,X_r^{s,x,\varepsilon})\)\dif r\\
&&-\sqrt{\varepsilon}\int_{s}^{t}\<\nabla\phi(X_r^{s,x,\varepsilon}),\sigma(r,X_r^{s,x,\varepsilon})\dif W_r\>,
\de
which together with $({\bf H}^{1'}_{b,\s})$ implies that
\ce
K^{s,x,\varepsilon}_t\leq |\phi(X_t^{s,x,\varepsilon})-\phi(x)|+CT+\left|\sqrt{\varepsilon}\int_{s}^{t}\<\nabla\phi(X_r^{s,x,\varepsilon}),\sigma(r,X_r^{s,x,\varepsilon})\dif W_r\>\right|.
\de
Furthermore, from the Burkholder-Davis-Gundy inequality and $({\bf H}^{1'}_{b,\s})$, it follows that
\ce
\mE\left[\sup\limits_{s \leq t \leq T}(K^{s,x,\varepsilon}_t)^p|\sF_s\right]&\leq& 3^{p-1}C\mE\left[\sup \limits_{s \leq t \leq T}|X_t^{s,x,\varepsilon}-x|^p|\sF_s\right]+3^{p-1}CT^p\no\\
&&+3^{p-1}C\mE\left[\left(\int_{s}^{T}|\nabla\phi(X_r^{s,x,\varepsilon})|^2\|\sigma(r,X_r^{s,x,\varepsilon})\|^2\dif r\right)^{p/2}|\sF_s\right]\no\\
&\leq&C.
\de	

By the same deduction to that for the above inequality, one can get (\ref{k2}).

Finally, in view of (\ref{k1}) and the Taylor formula, we obtain the required estimate (\ref{kexpo}). The proof is complete.
\end{proof}
	
\section{The LDP for generalized BSDEs}\label{bsdeLDP}

In this section, we prove the LDP for generalized BSDEs.

Consider Eq.(\ref{ygbsde}), i.e.
\be
\left\{\begin{array}{l}
\dif Y_t^{s,x,\varepsilon}=-f(t,X_t^{s,x,\varepsilon},Y_t^{s,x,\varepsilon},Z_t^{s,x,\varepsilon})\dif t-g(t,X_t^{s,x,\varepsilon},Y_t^{s,x,\varepsilon})\dif K_t^{s,x,\varepsilon}+Z_t^{s,x,\varepsilon}dW_t\\
Y_T^{s,x,\varepsilon}=h(X_T^{s,x,\varepsilon}),
\end{array}
\right.
\label{gbsde1}
\ee
where $h: \bar{\varTheta} \mapsto \mR^k$, $f: [0,T]\times\bar{\varTheta} \times\mR^k\times\mR^{k\times m} \mapsto \mR^k$, $g: [0,T]\times\bar{\varTheta} \times\mR^k \mapsto \mR^k$ are Borel measurable and $g$ belongs to $C^{1,2,2}([0,T]\times\bar{\varTheta} \times\mR^k, \mR^k)$.

Assume:
\begin{enumerate}[(${\bf H}_{f,g,h}$)]
\item $f(t,x,y,z)$ and $g(t,x,y)$ are continuous in $y$ and there exists a constant $L_3>0$ such that for $t\in [s,T]$, $x,x' \in \bar{\varTheta}$, $y,y' \in \mR^k$, $z,z' \in \mR^{k\times m}$,
\ce
&&|f(t,x,y,z)-f(t,x',y,z)|\leq L_3|x-x'|,\\
&&\<y-y', f(t,x,y,z)-f(t,x,y',z)\>\leq L_3|y-y'|^2,\\
&&|f(t,x,y,z)-f(t,x,y,z')|\leq L_3\|z-z'\|,\\
&&|g(t,x,y)-g(t,x',y)|\leq L_3|x-x'|,\\
&&\<y-y', g(t,x,y)-g(t,x,y')\>\leq L_3|y-y'|^2,\\
&&|h(x)-h(x')|\leq L_3|x-x'|,
\de
and
\ce
|f(t,x,y,z)|+|g(t,x,y)|\leq L_3(1+|y|+\|z\|). 
\de
\end{enumerate}

Under (${\bf H}_{f,g,h}$) and $({\bf H}^{1'}_{b,\s})$, by Theorem 1.6 in \cite{pz}, Eq.(\ref{gbsde1}) has a unique solution $\{Y_t^{s,x,\varepsilon},Z_t^{s,x,\varepsilon}\}_{0\leq s\leq t\leq T}$ satisfying
\ce
Y_t^{s,x,\varepsilon}&=&h(X_T^{s,x,\varepsilon})+\int_t^T f(r,X_r^{s,x,\varepsilon},Y_r^{s,x,\varepsilon},Z_r^{s,x,\varepsilon})\dif r+\int_t^T g(r,X_r^{s,x,\varepsilon},Y_r^{s,x,\varepsilon})\dif K_r^{s,x,\varepsilon}\no\\
&&-\int_t^T Z_r^{s,x,\varepsilon}\dif W_r,
\de
and
\ce
&&\mE\left[\sup \limits_{s \leq t \leq T} |Y_t^{s,x,\varepsilon}|^2+\int_s^T\|Z_r^{s,x,\varepsilon}\|^2\dif r \bigg{|} \sF_s \right]<C(1+|x|^2),\\
&&\sup \limits_{s \leq t \leq T}|Y_t^{s,x,\varepsilon}|\leq C, a.s.
\de
where the constant $C>0$ is independent of $x,\e$. Set 
\ce
u^\varepsilon(s,x):=Y_s^{s,x,\varepsilon}, \quad s\in[0,T],x\in \bar{\varTheta}
\de
and $u^\varepsilon(s,x)$ is deterministic since $Y_s^{s,x,\varepsilon}$ is $\mathscr{F}^s_s$-measurable and $\mathscr{F}^s_s=\{\emptyset,\Omega\}$. Then by Theorem 5.1 in \cite{pr1}, we know that $u^\e(s,x)=Y_s^{s,x,\varepsilon}$ is continuous in $(s,x)$. Moreover, the Markov property of the solution $X^{s,x,\varepsilon}$ for Eq.(\ref{sdef1}) and the uniqueness of the solution $Y^{s,x,\varepsilon}$ for Eq.(\ref{gbsde1}) imply that for $s\leq t\leq T$
\ce
u^\varepsilon(t,X_t^{s,x,\varepsilon})=Y_t^{t,X_t^{s,x,\varepsilon},\varepsilon}=Y_t^{s,x,\varepsilon}.
\de

Next, for any $\varrho \in C([0,T],\bar{\varTheta})$, set 
$$
\Pi^\varepsilon(\varrho)(\cdot):=u^\varepsilon(\cdot,\varrho_{\cdot}),
$$
and it holds that $\Pi^\e$ is a mapping from $C([0,T], \bar{\varTheta})$ to $C([0,T], \bar{\varTheta})$ and $Y^{s,x,\varepsilon}=\Pi^\varepsilon(X^{s,x,\varepsilon})$. Moreover, we have the following result.
				
\bl\label{piecont}
Under $({\bf H}_{f,g,h})$ and $({\bf H}^{1'}_{b,\s})$, for any $\e>0$, $\Pi^\varepsilon$ is continuous.
\el
\begin{proof}
First of all, we take $\{\varrho^n, \varrho: n\in\mN\}\subset C([s,T], \bar{\varTheta})$ such that $\varrho^n$ converges to $\varrho$. Then there exists a constant $M>0$ such that for all $n\in \mN$,
\ce
\|\varrho^n\|_\infty\leq M, \quad \|\varrho\|_\infty\leq M.
\de
Since $u^\varepsilon(t,x)$ is continuous in $(t,x)$, $u^\varepsilon$ is uniformly continuous on any compact subset of $[s,T]\times (B(0,M)\cap\bar{\varTheta})$, where $B(0,M):=\{x\in\mR^d: |x|\leq M\}$. Therefore, for any $\eta>0$, there exists $\delta>0$ such that for $v, v'\in[s,T], |v-v'|<\delta$ and $y,y'\in B(0,M)\cap\bar{\varTheta}, |y-y'|<\delta$, 
$$
|u^\varepsilon(v,y)-u^\varepsilon(v',y')|\leq\eta.
$$

Besides, since $\varrho^n$ converges to $\varrho$, there exists a $N>0$ such that for all $n\geq N$, $\|\varrho^n-\varrho\|_\infty \leq \delta$. 

Collecting the above deduction, we obtain that for all $n\geq N$, 
$$
\|\Pi^\e(\varrho^n)-\Pi^\e(\varrho)\|_{\infty}=\sup \limits_{v\in [s,T]}|u^\varepsilon(v,\varrho^n_v)-u^\varepsilon(v,\varrho_v)|\leq\eta.
$$		
The proof is complete.
\end{proof}
	
In order to study the LDP about $Y^{s,x,\varepsilon}$, we introduce the following auxiliary backward differential equation:
\be
\left\{\begin{array}{l}
\dif \psi_t^{s,x}=-f(t,\varphi_t^{s,x},\psi_t^{s,x},0)\dif t-g(t,\varphi_t^{s,x},\psi_t^{s,x})\dif K_t^{s,x},\\
\psi_T^{s,x}=h(\varphi_T^{s,x}).
\end{array}
\right.
\label{bde1}
\ee
Under $({\bf H}_{f,g,h})$ and $({\bf H}^{1'}_{b,\s})$, the above equation has a unique solutions $\psi^{s,x}$ satisfying
\ce
\psi_t^{s,x}=h(\varphi_T^{s,x})+\int_t^T f(r,\varphi_r^{s,x},\psi_r^{s,x},0)dr+\int_t^T g(r,\varphi_r^{s,x},\psi_r^{s,x})\dif K_r^{s,x},
\de
and $\sup \limits_{s \leq t \leq T}|\varphi_t^{s,x}|\leq C$. Then, set
$$
u(s,x):=\psi_s^{s,x},\quad s\in[0,T],x\in \bar{\varTheta},
$$
and $u$ is well defined on $[0,T]\times\bar{\varTheta}$. Moreover, by the uniqueness of the solutions for Eq.(\ref{def1}) and Eq.(\ref{bde1}), it holds that $u(t,\varphi_t^{s,x})=\psi_t^{t,\varphi_t^{s,x}}=\psi_t^{s,x}$. 

In the following, we define the mapping $\Pi$ by
\ce
\Pi(\varrho)(\cdot):=u(\cdot,\varrho_\cdot) \quad for\enspace all\enspace \varrho \in C([0,T];\bar{\varTheta})
\de
and $\Pi$ is well defined. In order to describe the relation between $\Pi^\e$ and $\Pi$, we prepare the following estimate.
	
\bl\label{yshou}
Under $({\bf H}_{f,g,h})$ and $({\bf H}^{1'}_{b,\s})$,  for all $\varepsilon \in (0,1)$, there exists a constant $C>0$ independent of $x$ and $\varepsilon$, such that
\ce
\sup\limits_{s\leq t\leq T}\mE\big[|Y_t^{s,x,\varepsilon}-\psi_t^{s,x}|^4\big] \leq C\varepsilon.
\de
\el
\begin{proof}
For any $\lambda>0$ and $s\leq t\leq T$, applying It\^o's formula to $e^{\lambda K^{s,x,\varepsilon}_t}|Y_t^{s,x,\varepsilon}-\psi_t^{s,x}|^2$, we get 
\ce
&&e^{\lambda K^{s,x,\varepsilon}_t}|Y_t^{s,x,\varepsilon}-\psi_t^{s,x}|^2+\int_{t}^{T} \lambda e^{\lambda K^{s,x,\varepsilon}_r} |Y_r^{s,x,\varepsilon}-\psi_r^{s,x}|^2\dif K^{s,x,\varepsilon}_r+\int_t^T  e^{\lambda K^{s,x,\varepsilon}_r}\|Z_r^{s,x,\varepsilon}\|^2 \dif r\no\\
&=& e^{\lambda K^{s,x,\varepsilon}_T}|h(X_T^{s,x,\varepsilon})-h(\varphi_T^{s,x})|^2- 2\int_t^T  e^{\lambda K^{s,x,\varepsilon}_r}\<Y_r^{s,x,\varepsilon}-\psi_r^{s,x}, Z_r^{s,x,\varepsilon}\dif W_r\>\no\\
&&+ 2\int_t^T e^{\lambda K^{s,x,\varepsilon}_r}\<Y_r^{s,x,\varepsilon}-\psi_r^{s,x},f(r,X_r^{s,x,\varepsilon},Y_r^{s,x,\varepsilon},Z_r^{s,x,\varepsilon})-f(r,\varphi_r^{s,x},\psi_r^{s,x},0)\>\dif r\no\\
&&+ 2\int_t^T e^{\lambda K^{s,x,\varepsilon}_r}\<Y_r^{s,x,\varepsilon}-\psi_r^{s,x},g(r,X_r^{s,x,\varepsilon},Y_r^{s,x,\varepsilon})-g(r,\varphi_r^{s,x},\psi_r^{s,x})\>\dif K^{s,x,\varepsilon}_r\\
&&+2\int_t^T e^{\lambda K^{s,x,\varepsilon}_r}\<Y_r^{s,x,\varepsilon}-\psi_r^{s,x}, g(r,\varphi_r^{s,x},\psi_r^{s,x})\>\dif (K^{s,x,\varepsilon}_r-K^{s,x}_r).
\de
Then by taking a suitable $\l$, (${\bf H}_{f,g,h}$) and the elementary inequality imply that
\ce
&&e^{\lambda K^{s,x,\varepsilon}_t}|Y_t^{s,x,\varepsilon}-\psi_t^{s,x}|^2+\frac{1}{2}\int_t^T  e^{\lambda K^{s,x,\varepsilon}_r}\|Z_r^{s,x,\varepsilon}\|^2 \dif r\\
&\leq& L^2_3e^{\lambda K^{s,x,\varepsilon}_T}|X_T^{s,x,\varepsilon}-\varphi_T^{s,x}|^2- 2\int_t^T  e^{\lambda K^{s,x,\varepsilon}_r}\<Y_r^{s,x,\varepsilon}-\psi_r^{s,x},Z_r^{s,x,\varepsilon}\dif W_r\>\no\\
&&+C\int_t^T e^{\lambda K^{s,x,\varepsilon}_r}|Y_r^{s,x,\varepsilon}-\psi_r^{s,x}|^2\dif r+\int_t^T e^{\lambda K^{s,x,\varepsilon}_r}|X_r^{s,x,\varepsilon}-\varphi_r^{s,x}|^2\dif r\\
&&+\int_t^T e^{\lambda K^{s,x,\varepsilon}_r}|X_r^{s,x,\varepsilon}-\varphi_r^{s,x}|^2\dif K^{s,x,\varepsilon}_r\\
&&+2\int_t^T e^{\lambda K^{s,x,\varepsilon}_r}\<Y_r^{s,x,\varepsilon}-\psi_r^{s,x},g(r,\varphi_r^{s,x},\psi_r^{s,x})\>\dif (K^{s,x,\varepsilon}_r-K^{s,x}_r),
\de
which yields that
\be
&&|Y_t^{s,x,\varepsilon}-\psi_t^{s,x}|^2+\frac{1}{2}\mE\left[\int_t^T  e^{\lambda K^{s,x,\varepsilon}_r}\|Z_r^{s,x,\varepsilon}\|^2 \dif r\bigg{|}\sF_t\right]\no\\
&\leq&\mE\left[L^2_3e^{\lambda K^{s,x,\varepsilon}_T}|X_T^{s,x,\varepsilon}-\varphi_T^{s,x}|^2\bigg{|}\sF_t\right]+\mE\left[\int_t^T e^{\lambda K^{s,x,\varepsilon}_r}|X_r^{s,x,\varepsilon}-\varphi_r^{s,x}|^2\dif r\bigg{|}\sF_t\right]\no\\
&&+\mE\left[\int_t^T e^{\lambda K^{s,x,\varepsilon}_r}|X_r^{s,x,\varepsilon}-\varphi_r^{s,x}|^2\dif K^{s,x,\varepsilon}_r\bigg{|}\sF_t\right]+C\mE\left[\int_t^T e^{\lambda K^{s,x,\varepsilon}_r}|Y_r^{s,x,\varepsilon}-\psi_r^{s,x}|^2\dif r\bigg{|}\sF_t\right]\no\\
&&+2\mE\left[\int_t^T e^{\lambda K^{s,x,\varepsilon}_r}\<Y_r^{s,x,\varepsilon}-\psi_r^{s,x}, g(r,\varphi_r^{s,x},\psi_r^{s,x})\>\dif (K^{s,x,\varepsilon}_r-K^{s,x}_r)\bigg{|}\sF_t\right].
\label{yzesti}
\ee

Next, we estimate the last term in the above inequality. By It\^o's formula to $e^{\lambda K^{s,x,\varepsilon}_t}\<Y_t^{s,x,\varepsilon}-\psi_t^{s,x},g(t,\varphi_t^{s,x},\psi_t^{s,x})\>(K^{s,x,\varepsilon}_t-K^{s,x}_t)$, it holds that
\ce
&& \int_t^T e^{\lambda K^{s,x,\varepsilon}_r}\<Y_r^{s,x,\varepsilon}-\psi_r^{s,x}, g(r,\varphi_r^{s,x},\psi_r^{s,x})\>\dif (K^{s,x,\varepsilon}_r-K^{s,x}_r)\\
&=& e^{\lambda K^{s,x,\varepsilon}_T}\<Y_T^{s,x,\varepsilon}-\psi_T^{s,x}, g(T,\varphi_T^{s,x},\psi_T^{s,x})\>(K^{s,x,\varepsilon}_T-K^{s,x}_T)\no\\
&&-e^{\lambda K^{s,x,\varepsilon}_t}\<Y_t^{s,x,\varepsilon}-\psi_t^{s,x}, g(t,\varphi_t^{s,x},\psi_t^{s,x})\>(K^{s,x,\varepsilon}_t-K^{s,x}_t)\no\\
&&+ \int_t^T \Bigg\{\<f(r,X_r^{s,x,\varepsilon},Y_r^{s,x,\varepsilon},Z_r^{s,x,\varepsilon})-f(r,\varphi_r^{s,x},\psi_r^{s,x},0), g(r,\varphi_r^{s,x},\psi_r^{s,x})\>\no\\
&&\quad-\Bigg<Y_r^{s,x,\varepsilon}-\psi_r^{s,x}, \big[\frac{\partial}{\partial{r}}g(r,\varphi_r^{s,x},\psi_r^{s,x})+\nabla_x g(r,\varphi_r^{s,x},\psi_r^{s,x})\cdot b(r,\varphi_r^{s,x})\no\\
&&\quad\quad-\nabla_y g(r,\varphi_r^{s,x},\psi_r^{s,x})\cdot f(r,\varphi_r^{s,x},\psi_r^{s,x},0)\big]\Bigg>\Bigg\}e^{\lambda K^{s,x,\varepsilon}_r}(K^{s,x,\varepsilon}_r-K^{s,x}_r)\dif r\no\\
&&- \int_t^T e^{\lambda K^{s,x,\varepsilon}_r}(K^{s,x,\varepsilon}_r-K^{s,x}_r)\<g(r,\varphi_r^{s,x},\psi_r^{s,x}), Z_r^{s,x,\varepsilon}\dif W_r\>\no\\
&&+ \int_t^T \big[\<g(r,X_r^{s,x,\varepsilon},Y_r^{s,x,\varepsilon}), g(r,\varphi_r^{s,x},\psi_r^{s,x})\>-\lambda\<Y_r^{s,x,\varepsilon}-\psi_r^{s,x}, g(r,\varphi_r^{s,x},\psi_r^{s,x})\>\big]\cdot\no\\
&&\quad e^{\lambda K^{s,x,\varepsilon}_r}(K^{s,x,\varepsilon}_r-K^{s,x}_r)\dif K^{s,x,\varepsilon}_r\no\\
&&+ \int_t^T \Big[-|g(r,\varphi_r^{s,x},\psi_r^{s,x})|^2-\<Y_r^{s,x,\varepsilon}-\psi_r^{s,x}, \big(\nabla_x g(r,\varphi_r^{s,x},\psi_r^{s,x})\nabla\phi(\varphi_r^{s,x})\no\\
&&\quad-\nabla_y g(r,\varphi_r^{s,x},\psi_r^{s,x})\cdot g(r,\varphi_r^{s,x},\psi_r^{s,x})\big)\>\Big] e^{\lambda K^{s,x,\varepsilon}_r}(K^{s,x,\varepsilon}_r-K^{s,x}_r)\dif K^{s,x}_r.
\de 
Moreover, by (${\bf H}_{f,g,h}$), we have that
\be
&&2\mE \left[\int_t^T e^{\lambda K^{s,x,\varepsilon}_r}\<Y_r^{s,x,\varepsilon}-\psi_r^{s,x}, g(r,\varphi_r^{s,x},\psi_r^{s,x})\>\dif (K^{s,x,\varepsilon}_r-K^{s,x}_r) |\sF_t\right]\no\\
&\leq& C\mE \Big[e^{\lambda K^{s,x,\varepsilon}_T} |X_T^{s,x,\varepsilon}-\varphi_T^{s,x}|\cdot|K^{s,x,\varepsilon}_T-K^{s,x}_T||\sF_t\Big]\no\\
&&+e^{\lambda K^{s,x,\varepsilon}_t}|Y_t^{s,x,\varepsilon}-\psi_t^{s,x}|\cdot|K^{s,x,\varepsilon}_t-K^{s,x}_t| \no\\
&&+ C\mE \Big[\int_t^T \big( |X_r^{s,x,\varepsilon}-\varphi_r^{s,x}|+|Y_r^{s,x,\varepsilon}-\psi_r^{s,x}|+\|Z_r^{s,x,\varepsilon}\|\big) e^{\lambda K^{s,x,\varepsilon}_r}|K^{s,x,\varepsilon}_r-K^{s,x}_r|\dif r|\sF_t\Big]\no\\
&&+C\mE \Big[\int_t^T \big(1+|Y_r^{s,x,\varepsilon}-\psi_r^{s,x}|\big) e^{\lambda K^{s,x,\varepsilon}_r}|K^{s,x,\varepsilon}_r-K^{s,x}_r|\dif K^{s,x,\varepsilon}_r|\sF_t\Big]\no\\
&&+C\mE \Big[\int_t^T \big(1+|Y_r^{s,x,\varepsilon}-\psi_r^{s,x}|\big) e^{\lambda K^{s,x,\varepsilon}_r}|K^{s,x,\varepsilon}_r-K^{s,x}_r|\dif K^{s,x}_r|\sF_t\Big]\no\\
&\leq&C\mE \Big[|X_T^{s,x,\varepsilon}-\varphi_T^{s,x}|^2|\sF_t\Big]+ C\mE \Big[ e^{2\lambda K^{s,x,\varepsilon}_T}|K^{s,x,\varepsilon}_T-K^{s,x}_T|^2|\sF_t\Big]\no\\
&&+\frac{1}{2} |Y_t^{s,x,\varepsilon}-\psi_t^{s,x}|^2+Ce^{2\lambda K^{s,x,\varepsilon}_t}|K^{s,x,\varepsilon}_t-K^{s,x}_t|^2\no\\
&&+C\mE \Big[\int_t^T \big( |X_r^{s,x,\varepsilon}-\varphi_r^{s,x}|^2+|Y_r^{s,x,\varepsilon}-\psi_r^{s,x}|^2\big) \dif r|\sF_t\Big]\no\\
&&+C\mE \Big[\int_t^T e^{2\lambda K^{s,x,\varepsilon}_r}|K^{s,x,\varepsilon}_r-K^{s,x}_r|^2\dif r|\sF_t\Big]\no\\
&&+\frac{1}{4}\mE \Big[\int_t^T e^{\lambda K^{s,x,\varepsilon}_r} \|Z_r^{s,x,\varepsilon}\|^2\dif r|\sF_t\Big]+C\mE \Big[\int_t^T e^{\lambda K^{s,x,\varepsilon}_r}|K^{s,x,\varepsilon}_r-K^{s,x}_r|^2\dif r|\sF_t\Big]\no\\
&&+C\mE \Big[e^{\lambda K^{s,x,\varepsilon}_T}K^{s,x,\varepsilon}_T\sup\limits_{s\leq r\leq T}|K^{s,x,\varepsilon}_r-K^{s,x}_r| |\sF_t\Big]\no\\
&&+C\mE \Big[e^{\lambda K^{s,x,\varepsilon}_T}K^{s,x}_T\sup\limits_{s\leq r\leq T}|K^{s,x,\varepsilon}_r-K^{s,x}_r| |\sF_t\Big].
\label{gesti}
\ee

Combining (\ref{yzesti}) with (\ref{gesti}), by the Jensen inequality, the H\"older inequality, Lemma \ref{xshou} and \ref{kshou}, (\ref{k1}), (\ref{k2}) and (\ref{kexpo}), we obtain that
\ce
\mE \big[|Y_t^{s,x,\varepsilon}-\psi_t^{s,x}|^4|\sF_s\big] \leq C\varepsilon+C\mE \big[\int_t^T |Y_r^{s,x,\varepsilon}-\psi_r^{s,x}|^4\dif r|\sF_s\big].
\de
Consequently, it follows from Gronwall's lemma that
\ce
\mE \big[|Y_t^{s,x,\varepsilon}-\psi_t^{s,x}|^4|\sF_s\big] \leq C\varepsilon.
\de	
Then, we have that for all $t\in[s,T]$, 
\ce
\mE \big[|Y_t^{s,x,\varepsilon}-\psi_t^{s,x}|^4\big] &=& \mE\Big[\mE \big[|Y_t^{s,x,\varepsilon}-\psi_t^{s,x}|^4|\sF_s\big]\Big]\no\\ 
&\leq& C\varepsilon.
\de	
The proof is complete.
\end{proof}

\br\label{limi}
Note that $u^\varepsilon(s,x):=Y_s^{s,x,\varepsilon}$ is a unique viscosity solution of PDE (\ref{pde1}), i.e.
\ce\left\{\begin{array}{l}
\frac{\p u_i^\e(s,x)}{\p s}+b_j(s,x)\p_j u_i^\e(s,x)+\frac{\e}{2}(\s\s^*)_{jl}(s,x) \p_j\p_l u_i^\e(s,x)+f_i(s, x, u^\e(s,x), (\triangledown u_i^\e \s)(s,x) )=0,\\
 ~(s,x)\in[0,T]\times\bar{\varTheta}, \quad i=1, \cdots, k\\ 
\frac{\p u_i^\e(s,x)}{\p n}+g_i(s,x,u^\e(s,x))=0, \\
~(s,x)\in[0,T]\times\p \bar{\varTheta}, \quad i=1, \cdots, k\\
u^\e(T,x)=h(x), \quad x\in\bar{\varTheta}.
\end{array}
\right.
\de
By the above lemma, we could think that as $\e$ tends to $0$, $u^\varepsilon$ converges to $u$ which solves the following PDE:
\ce\left\{\begin{array}{l}
\frac{\p u_i(s,x)}{\p s}+b_j(s,x)\p_j u_i(s,x)+f_i(s, x, u(s,x), 0)=0,\\
 ~(s,x)\in[0,T]\times\bar{\varTheta}, \quad i=1, \cdots, k\\ 
\frac{\p u_i(s,x)}{\p n}+g_i(s,x,u(s,x))=0, \\
~(s,x)\in[0,T]\times\p \bar{\varTheta}, \quad i=1, \cdots, k\\
u(T,x)=h(x), \quad x\in\bar{\varTheta}.
\end{array}
\right.
\de
\er

The following lemma characterizes the relationship between $\Pi^\e$ and $\Pi$.
	
\bl\label{ydpl2}
Under $({\bf H}_{f,g,h})$ and $({\bf H}^{1'}_{b,\s})$, the maps $\Pi^\varepsilon$ converge uniformly to $\Pi$ on every compact subset of $C([s,T];\bar{\varTheta})$.
\el
\begin{proof}
Let $\Upsilon$ be any compact subset of $C([s,T];\bar{\varTheta})$. So, for $\varrho\in \Upsilon$, there exists a constant $M>0$ such that $\sup\limits_{\varrho \in \Upsilon}\sup \limits_{v\in [s,T]} |\varrho_v|\leq M$. Moreover, it holds that
\ce
\|\Pi^\varepsilon(\varrho)-\Pi(\varrho)\|_\infty &=& \sup \limits_{v\in [s,T]} |\Pi^\varepsilon(\varrho)(v)-\Pi(\varrho)(v)|\no\\
&=& \sup \limits_{v\in [s,T]} |u^\varepsilon(v,\varrho_v)-u(v,\varrho_v)|\no\\
&=& \sup \limits_{v\in [s,T]} |Y_v^{v,\varrho_v,\varepsilon}-\psi_v^{v,\varrho_v}|,
\de
and
\ce
\sup \limits_{\varrho \in \Upsilon}\lVert \Pi^\varepsilon(\varrho)-\Pi(\varrho) \rVert_\infty^4
&=& \sup\limits_{\varrho \in \Upsilon}\sup \limits_{v\in [s,T]} |Y_v^{v,\varrho_v,\varepsilon}-\psi_v^{v,\varrho_v}|^4\no\\
&\leq& \sup\limits_{|x| \leq M}\sup \limits_{v\in [s,T]}|Y_v^{v,x,\varepsilon}-\psi_v^{v,x}|^4,
\de 
which together with Lemma \ref{yshou} yields that		
\ce 
\sup \limits_{\varrho \in \Upsilon}\lVert \Pi^\varepsilon(\varrho)-\Pi(\varrho) \rVert_\infty^4 &\leq& \sup\limits_{|x| \leq M}\sup \limits_{v\in [s,T]}\mE|Y_v^{v,x,\varepsilon}-\psi_v^{v,x}|^4\no\\
&\leq& \sup\limits_{|x| \leq M}\sup \limits_{v\in [s,T]}\sup\limits_{v \leq t \leq T} \mE \big[|Y_t^{v,x,\varepsilon}-\psi_t^{v,x}|^4 \big]\no\\
&\leq& C\varepsilon.
\de
The proof is complete.		
\end{proof}

Now, we state and prove the main result in this section.

\bt\label{yldp} 
Assume that $({\bf H}_{f,g,h})$, $({\bf H}^{1'}_{b,\s})$ and $({\bf H}^2_{\s})$ hold. Then $\{Y^{s,x,\varepsilon}, \e>0\}$ satisfies, as $\varepsilon$ goes to 0, a large deviation principle with the following good rate function 
\ce
S'(\Gamma)=\inf\{S(\Psi): \Gamma_t=\Pi(\Psi)(t)=u(t,\Psi_t) ~\mbox{for all}~t\in [s,T] \}.
\de
\et
\begin{proof}
By Theorem $\ref{fldp}$, we know that $\{X^{s,x,\varepsilon}, \e>0\}$ satisfies the LDP with the good rate function 
\ce
S(\Psi)
=\left\{\begin{array}{ll}
\frac{1}{2} \inf\limits_{\Phi \in F^{-1}(\Psi)} \int_s^T (\dot{\Phi}_t-b(t,\Psi_t))^* (\s\s^*)^{-1}(t,\Psi_t) (\dot{\Phi}_t-b(t,\Psi_t)) \dif t, \\
\qquad\quad\text{if $\Phi_t$ is absolutely continuous, $\Phi_s=x$ and $F^{-1}(\Psi) \neq \emptyset$};\\
+\infty, \quad \text{in the opposite case}. 
\end{array}
\right. 
\de

Besides, we notice $Y^{s,x,\varepsilon}=\Pi^\varepsilon(X^{s,x,\varepsilon})$. And by Lemma $\ref{piecont}$ and $\ref{ydpl2}$, one can obtain that $\Pi^\varepsilon$ is continuous and converges uniformly to $\Pi$ on every compact subset of $C([s,T];\bar{\varTheta})$. Therefore, from Theorem \ref{yasuo}, it follows that the process $Y^{s,x,\varepsilon}$ satisfies, as $\varepsilon$ goes to 0, a LDP with the following good rate function 
\ce
S'(\Gamma)=\inf\{S(\Psi): \Gamma_t=\Pi(\Psi)(t)=u(t,\Psi_t)\; for\, all\; t\in [s,T] \}.
\de
The proof is complete.
\end{proof}

\end{document}